\documentclass{amsart}

\newtheorem{theorem}{Theorem}[section]
\newtheorem{lemma}[theorem]{Lemma}
\newtheorem{corollary}[theorem]{Corollary}

\theoremstyle{definition}

\theoremstyle{remark}

\newtheorem{conjecture}[theorem]{Conjecture}
\newtheorem{setup}[theorem]{Setup}

\numberwithin{equation}{section}

\begin{document}

\title[On free subalgebras of varieties]{On free subalgebras of varieties}

\author[Renato Fehlberg J\'unior]{Renato Fehlberg J\'unior}\thanks{The authors would like to thank Universidade Federal do Esp\'{\i}rito Santo and Universidade de S\~ao Paulo for the hospitality and support during respective short term visits. Also, the authors would like to thank professor Ivan Shestakov for valuable comments.}
\address{Department of Mathematics, Universidade Federal do Esp\'{\i}rito Santo, Vit\'{o}ria, ES, Brazil}
\email{renato.fehlberg@ufes.br}

\author[J. S\'anchez]{Javier S\'anchez}\thanks{Partially supported by grant \#2015/09162-9,
S\~ao Paulo Research Foundation (FAPESP) - Brazil, and CNPq-Brazil Grant 305773/2018-6 .}
\address{Department of Mathematics, University of S\~{a}o Paulo, S\~{a}o Paulo, SP,
05508-090, Brazil}
\email{jsanchez@ime.usp.br}


\date{\today}


\begin{abstract}
We show that some results of L. Makar-Limanov, P. Malcolmson and Z. Reichstein on the existence of free associative algebras are valid in the more general context of varieties of algebras.
\end{abstract}

\maketitle

\section*{Introduction}

As stated in \cite[Conjecture~1.1]{Agata}, L. Makar-Limanov made the following conjecture:

\begin{conjecture}
Let $K$ be a field, $A$ be an associative $K$-algebra and $F$ be a field extension of $K$.
	If $F\otimes_KA$ contains a free $K$-algebra on at least two free generators, then $A$ also contains
	a free $K$-algebra on the same number of free generators.
\end{conjecture}

In \cite[Theorem~1(b)]{ZI}, Z. Reichstein proved that Makar-Limanov's conjecture holds true when
the field $K$ is uncountable:
\begin{theorem}\label{theo:Reichstein}
Let $K$ be an uncountable field, $A$ an associative $K$-algebra and $F$ a field extension of $K$. If 
$F\otimes_KA$ contains a copy of a free (noncommutative) associative $K$-algebra, then so does $A$.
\end{theorem}
 In his proof, Z. Reichstein made essential use of the following result by
L. Makar-Limanov and P. Malcolmson \cite[Lemma~1]{MLM}:

\begin{lemma}\label{lem:MLM}
Suppose that $K$ is a field
with prime subfield $K_0$ and $A$ is an associative $K$-algebra. Then $x_1,\dotsc,x_n\in A$
are the free generators of a noncommutative free $K$-subalgebra if, and only if, they are the free generators of a
free $K_0$-subalgebra.
\end{lemma}

On the other hand, A. Smoktunowicz proved in \cite[Theorem~1.2]{Agata} that this conjecture fails when $K$ is a countable
field. More precisely, she showed that for every countable field $K$, there is
an associative $K$-algebra $A$ without free subalgebras on at least two free generators and a
field extension $F$ of $K$ such that the algebra $F\otimes_K A$ contains a  free $K$-algebra
on at least two free generators.

 \medskip
 
The main aim of this paper is to illustrate the fact that  similar phenomena about the existence of free algebras hold true in the context of varieties of (not necessarily associative) algebras.

Before giving more details, we fix some notation that will be used throughout. For unexplained terminology, the reader is referred to \cite[Chapter~1]{Ringsthatarenearly}.

Let $Y$ be a set and $K$ be a field. Let $X$ be a countable set of symbols 
$X=\{x_1,x_2,\dotsc\}$ and
let $\mathfrak{M}$ be a variety of $K$-algebras with defining identities $I\subseteq K\{X\}$.

By $K\{Y\}$, we denote the free (nonassociative) $K$-algebra on $Y$.  Thus, for any $K$-algebra $A$ and map $\theta\colon Y\rightarrow A$ there exists a unique
homomorphism $\Theta\colon K\{Y\}\rightarrow A$ which extends $\theta$.

 We will denote by $K_\mathfrak{M}\{Y\}$ the free $K$-algebra in the variety $\mathfrak{M}$ with set of free generators $Y$. Thus, for any $K$-algebra $A\in\mathfrak{M}$ and map $\theta\colon Y\rightarrow A$ there exists a unique
 homomorphism of $K$-algebras $\Theta\colon K_\mathfrak{M}\{Y\}\rightarrow A$ which extends $\theta$.

 If 
 $K\subseteq F$ is a field extension, then   $\mathfrak{M}_F$ denotes the variety of $F$-algebras with defining identities $I$. Then $F\otimes_K K_\mathfrak{M}\{Y\}$ is the free $F$-algebra of $\mathfrak{M}_F$ with set of free generators
$\{1\otimes y_1,\dotsc,1\otimes y_n\}$, if $Y=\{y_1,\dotsc,y_n\}$.

We will only consider homogeneous varieties $\mathfrak{M}$ of $K$-algebras. Hence, if $A\in \mathfrak{M}$, then $K$-algebra
$F\otimes_K A\in \mathfrak{M}$.

\medskip


Let $K$ be a field and $\mathfrak{M}$ be a homogeneous variety of $K$-algebras.
We say that  $\mathfrak{M}$ is an \emph{MLM variety} if 
for any field extension $K\subseteq F$, any $A\in\mathfrak{M}_F$, and subset of at least two elements
$Y=\{y_1,\dotsc,y_n\}\subseteq A$  such
that the $K$-subalgebra of $A$ generated by $Y$ is the free $K$-algebra in the variety  $\mathfrak{M}$ with
set of free
generators $Y$, then the $F$-subalgebra of $A$ generated by $Y$ is the free $F$-algebra in the variety $\mathfrak{M}_F$
with set of free generators $Y$. The name MLM stands for
Makar-Limanov and Malcolmson. Note that Lemma~\ref{lem:MLM} states that the variety of associative $K$-algebras is MLM. 

Suppose that $K$ is an uncountable field and $\mathfrak{M}$ is a homogeneous variety of $K$-algebras. We say that   $\mathfrak{M}$ is a \emph{Reichstein variety} if for any $A\in \mathfrak{M}$ and field extension
$F$ of $K$   such that $F\otimes_K A$ contains a free $K$-algebra in the variety $\mathfrak{M}$ on at least
two free generators, then $A$ contains a free $K$-algebra in the variety $\mathfrak{M}$ on the same
number of free generators. Note that Theorem~\ref{theo:Reichstein} shows that the variety of associative $K$-algebras is Reichstein.

In Section~\ref{conjecture}, we prove that if $K$ is an uncountable field, MLM varieties of $K$-algebras are Reichstein. 
Also, as an easy consequence of \cite{Agata}, we show that if $K$ is a countable field and $\mathfrak{M}$
is either the variety of $K$-algebras generated by the special Jordan $K$-algebras or the
variety of Lie $K$-algebras, then  Smoktunowicz's result holds. That is, there exist
a $K$-algebra $A$ in $\mathfrak{M}$ and a field extension $F$ of $K$ such that $A$ does
not contain a free $K$-algebra in $\mathfrak{M}$ on at least two free generators but
$F\otimes_K A$ contains a free $K$-algebra in $\mathfrak{M}$ on at least two free generators.

In Section~\ref{ESFS}, we show that if $K$ is a field, the following homogeneous varieties of $K$-algebras are MLM:
\begin{itemize}
\item The variety $\mathcal{A}_K$ of all $K$-algebras
\item The variety $\mathcal{L}_K$ of all Lie $K$-algebras
\item The variety $\mathcal{C}_K$ of all commutative $K$-algebras
\item The variety $\mathcal{AC}_K$ of all anticommutative $K$-algebras
\item The variety $\mathcal{SJ}_K$ generated by the special Jordan $K$-algebras.
\end{itemize}

We end this introduction showing that not all homogeneous varieties of $K$-algebras are MLM. For example, the variety $\mathcal{T}riv_K$
of $K$-algebras which satisfy the identity $x_1x_2=0$ is not MLM. 
This variety can be identified with the class of $K$-vector spaces, and the $K$-basis of these vector spaces
with free set of generators.
Let now $F$ be a nontrivial field extension of $K$,
and consider $A$, a one-dimensional $F$-algebra with basis $\{z\}\subset A$. Note that if $f_1,f_2\in F$ are $K$-linearly independent, then the $K$-subalgebra of $A$ generated by $Y=\{f_1z,f_2z\}$ is free on $Y$ of rank two. On the other hand
the $F$-algebra generated by $Y$ is not free on $Y$ because $f_1z,f_2z$ are $F$-linearly dependent. 
Another example of a variety that is not MLM is the variety of
commutative and associative $K$-algebras. Indeed, consider the field of fractions $F$ of the polynomial ring in two variables $K[x,y]$. Clearly, $F$
contains a free $K$-algebra on two generators, but $F$ does not contain a free $F$-algebra on $\{x,y\}$.



\section{MLM varieties are Reichstein}\label{conjecture}

In the first part of this
section we prove results analogous to the ones in \cite{ZI} in the context of varieties of algebras.
The proofs are natural adaptations of the ones by Z. Reichstein.

The proof of the following result can be found in \cite[Lemma~1]{ZI}.

\begin{lemma}\label{L3}
  Let $K$ be an uncountable field and let $X_1, X_2, \ldots $ be a countable number of Zariski closed subsets of $K^n$. If $\cup^\infty_{i=1}X_i=K^n$, then $X_i=K^n$ for some $i\geq 1$. \qed
\end{lemma}

Let $K$ be a field and $F$ be a field extension of $K$. Suppose that $\mathfrak{M}$
is a homogeneous variety of $K$-algebras.

Let $A\in\mathfrak{M}$.
If $z\in F$ and $a\in A$, we shall denote $za\in F\otimes_K A$ instead of $z\otimes a$.
Let $(a_{11},\dotsc,a_{1r_1})\in A^{r_1},\dotsc,(a_{n1},\dotsc,a_{nr_n})\in A^{r_n}$.
For $z_1=(z_{11},\dotsc, z_{1r_1})\in F^{r_1},\dotsc,\linebreak z_n=(z_{n1},\dotsc,z_{nr_n})\in F^{r_n}$, set
\begin{equation*}
    a_{z_i}=\sum_{j=1}^{r_i}z_{ij}a_{ij}\in F\otimes_K A, \quad i\in\{1,\dotsc,n\}.
\end{equation*}

\begin{lemma}\label{lem:Zariskiclosed}
Let $Y=\{y_1,\dotsc,y_n\}, n\geq 2$, be a finite set. Let $f_1,\dotsc,f_m\in K\{Y\}$ be polynomials in $n$ variables.
Then the $n$-tuples $(z_1,\dotsc,z_n)\in F^{r_1+\dotsb+r_n}$ such that
$f_1(a_{z_1},\dotsc,a_{z_n}),\dotsc, f_{m}(a_{z_1},\dotsc,a_{z_n})\in F\otimes_K A$ are $F$-linearly dependent, form
a Zariski closed subset of $F^{r_1+\dotsb+r_n}$ defined over $K$.
\end{lemma}

\begin{proof}
Let $d$ be the maximum of the degrees of $f_1,\dotsc,f_m$ and let
$e_1,e_2,\dotsc,e_s$ be a basis of the $K$-vector subspace of $A$ spanned by all the
possible evaluations in $\{a_{ij}\}_{i,j}$
of the monomials in $K\{Y\}$
 of degree $\leq d$. Notice that $e_1,\dotsc,e_s$ are also $F$-linearly independent in $F\otimes_K A$.
 Then, for $k=1,\dotsc,m$, we can write
 \begin{equation*}
     f_k(a_{z_1},\dotsc,a_{z_n})=\sum_{t=1}^s p_{kt}(z_1,\dotsc,z_n)e_t,
 \end{equation*}
where each $p_{kt}$ is an associative and commutative polynomial in $r_1+\dotsb+r_n$ variables with coefficients in $K$.

If $m> s$, then the set consisting of the $(z_1,\dotsc,z_n)\in F^{r_1+\dotsb+r_n}$ such that \linebreak
$f_1(a_{z_1},\dotsc,a_{z_n}),\dotsc,f_m(a_{z_1},\dotsc,a_{z_n})$ are $F$-linearly dependent equals
$F^{r_1+\dotsb+r_n}$.

Suppose now that $m\leq s$. Then $f_1(a_{z_1},\dotsc,a_{z_n}),\dotsc,f_m(a_{z_1},\dotsc,a_{z_n})$ are $F$-linearly dependent
if, and only if, the $m\times s$ matrix $(p_{kt}(z_1,\dotsc,z_n))_{k,t}$ has rank $\leq m-1$. This is
equivalent to the vanishing of the $m\times m$ minors of this matrix. Each minor is a commutative and associative polynomial
in $z_{11},\dotsc,z_{1r_1},\dotsc,z_{n1},\dotsc,z_{nr_n}$ with coeficients in $K$.
\end{proof}

With this lemmas we can proof the main result of this section.

\begin{theorem}\label{theo:MLMReichstein}
Let $Q\subseteq K\subseteq F$ be field extensions with $K$ uncountable, 
 $\mathfrak{M}$ be an MLM variety of $Q$-algebras and $A\in \mathfrak{M}_K$.
 If $F\otimes_KA$ contains a free $Q$-algebra in the variety $\mathfrak{M}$ with a finite set of free generators greater or equal than two,
 then so does $A$.
As a consequence, if $\mathfrak{M}$ is an MLM 
variety of $K$-algebras, then $\mathfrak{M}$ is a Reichstein variety of $K$-algebras. 
\end{theorem}

\begin{proof}
Suppose that the elements 
\begin{equation*}
    a_{u_i}=\sum_{j=1}^{r_i}u_{ij}a_{ij}\in F\otimes_K A,\quad i=1,\dotsc,n,
\end{equation*}
are the free generators of a free $Q$-algebra in the variety $\mathfrak{M}$ for some $(u_{11},\dotsc,u_{1r_1})\in F^{r_1},\dotsc,
(u_{n1},\dotsc,u_{nr_n})\in F^{r_n}$ and $a_{ij}\in A$. Since $\mathfrak{M}$ is MLM, these elements also generate a free $F$-algebra in the variety $\mathfrak{M}_F$ with free set of
generators $\{a_{u_1},\dotsc,a_{u_n}\}$.

Let $\{y_1,\dotsc,y_n\}$ be a finite set. Consider the free $K$-algebra $K\{y_1,\dotsc,y_n\}$
 and the
free algebra $K_\mathfrak{M}\{y_1,\dotsc,y_n\}$ in the variety $\mathfrak{M}$ with
free set of generators $\{y_1,\dotsc,y_n\}$.
Consider the natural homomorphism of $K$-algebras
\begin{equation*}
    \Phi\colon K\{y_1,\dotsc,y_n\}\rightarrow K_\mathfrak{M}\{y_1,\dotsc,y_n\}, y_i\mapsto y_i.
\end{equation*}

For each $d\geq 1$, fix monomials $m_{d1},\dotsc,m_{dt_d}$ be monomials of degree $d$ in the free algebra
$K\{y_1,\dotsc,y_n\}$ such that $\bigcup\limits_{d\geq 1} \{\Phi(m_{d1}),\dotsc,\Phi(m_{dt_d})\}$ is a $K$-basis of
$K_{\mathfrak{M}}\{y_1,\dotsc,y_n\}$.

For $p\geq 1$, let $X_p\subseteq F^{r_1+\dotsb+r_n}$ be the set of all $n$-tuples
\begin{equation*}
    ((z_{11},\dotsc,z_{1r_1}),\dotsc,(z_{n1},\dotsc,z_{nr_n}))
\end{equation*}
such that
$m_{d{l}}(a_{z_1},\dotsc,a_{z_n})$, $l=1,\dotsc,t_d$, $d\leq p$, are $F$-linearly dependent. We will write $r_1+\dotsb+r_n:=r$.
By Lemma~\ref{lem:Zariskiclosed}, $X_p$ is a closed subset of $F^r$ defined over $K$.

In order to prove the existence of the free algebra $K_\mathfrak{M}\{y_1,\dotsc,y_n\}$ in $A$, we must 
show that all $m_{dj}(a_{z_1},\dotsc,a_{z_n})$ are $K$-linearly independent for some $(z_1,\dotsc,z_n)\in K^r$.
Assume the contrary: for every $(z_1,\dotsc,z_n)\in K^r$, the elements $a_{z_1},\dotsc,a_{z_n}$
are such that there exists $p\geq 1$ with
\begin{equation*}
    \sum_{d\leq p,1\leq l\leq t_d} \lambda_{d_l}m_{d_l}(a_{z_1},\dotsc,a_{z_n})=0,\quad \textrm{for some }\lambda_{d_l}\in K.
\end{equation*}
In other words, $(z_1,\dotsc,z_n)\in X_p$. Hence $K^r=\bigcup\limits_{p\geq 1}X_p(K)$,
where  $$X_p(K)=\{(z_1,\dotsc,z_n)\in K^r: m_{d_{l}}(a_{z_1},\dotsc,a_{z_n}),\, l=1,\dotsc,t_d,\, d\leq p, \textrm{ are } K\textrm{-l.d.}\}.$$
By Lemma~\ref{L3}, $X_p(K)=K^r$ for some integer $p$.
Note that $X_p(K)\subseteq X_p$. Now, since $K^r$ is dense in $F^r$, we get that
$X_p=F^r$. A contradiction because $a_{u_i}$ generate a free algebra
$F_{\mathfrak{M}_F}\{a_{u_1},\dotsc,a_{u_n}\}$.

Now, if $Q=K$, one obtains the  last assertion of the theorem. 
\end{proof}

\medskip

We end this section showing that Theorem~\ref{theo:Reichstein} does not hold for the variety of
Lie $K$-algebras and the variety generated by the special Jordan $K$-algebras when the field $K$ is countable.

Suppose that $K$ is any countable field. By \cite[Theorem~1.4]{Agata}, there exists
a field extension $F$ of $K$ and a nil associative $K$-algebra $A$ such that the
associative algebra $F\otimes_K A$ contains a non-commutative free $K$-algebra on
a  set of free generators of at least two elements.

Consider the special Jordan $K$-algebra $A^{(+)}$. It is known that $A^{(+)}$ is a
special Jordan nil $K$-algebra. Hence $A^{(+)}$ does not contain a copy of $\mathcal{SJ}_K(Y)$, the free special Jordan $K$-algebra with set of free generators $Y$, where
$Y$ possesses at least two elements. Now since, $F\otimes_K A$ contains a non-commutative free $K$-algebra on
a  set of free generators of at least two elements,
the special Jordan algebra $F\otimes_K A^{(+)}$ contains a copy of $\mathcal{SJ}_K(Y)$
where $Y$ has at least two elements.

Now consider the Lie $K$-algebra $A^{(-)}$. It is known that $A^{(-)}$ is a
Lie Engel $K$-algebra. Hence $A^{(-)}$ does not contains a copy of $\mathcal{L}_K(Y)$ where
$Y$ possesses at least two elements. Now, since $F\otimes_K A$ contains a non-commutative free $K$-algebra on
a  set of free generators of at least two elements,
the Lie algebra $F\otimes_K A^{(-)}$ contains a copy of $\mathcal{L}_K(Y)$, the free Lie $K$-algebra with set of free generators $Y$,
where $Y$ has at least two elements.


\section{Examples of MLM varieties}\label{ESFS}

Our aim in this section is to show that some important varieties of algebras are MLM. Hence, when the ground field is uncountable, these varieties are Reichstein by Theorem~\ref{theo:MLMReichstein}, 

To avoid  repetitive arguments, we  establish here a Setup that will be used at the beginning of the proofs as a standard text.

\begin{setup}\label{setup}
Let $K\subseteq F$ be a field extension and
$\mathfrak{A}$ be a variety of $K$-algebras. Let
$A$ be an $F$-algebra in $\mathfrak{A}_F$. Suppose that $Y=\{y_1,\dotsc,y_n\}\subseteq A$
is a set of at least two elements such that the $K$-subalgebra of $A$
generated by $Y$ is $K_{\mathfrak{A}}\{Y\}$, the free $K$-algebra in the variety  $\mathfrak{A}$ with free set of
generators $Y$.

Consider the homomorphism of  $F$-algebras
\begin{align*}
  \mu \colon F\otimes_K K_\mathfrak{A}\{Y\} &\longrightarrow A\\
  c\otimes p &\mapsto cp
\end{align*}
If we want to prove that $\mathfrak{A}$ is an MLM variety, we must show that $\mu$ is injective. We proceed as in the proof of \cite[Lemma~1]{MLM}.
Clearly $\mu(c\otimes p)=0$ if, and only if, $c$ is zero or $p$ is zero. Either way, $c\otimes p$ is zero.
Hence suppose that
\begin{align}\label{eq:muzero}
\mu\left(\sum_{i=1}^{n}{c_i\otimes p_i}\right)=0,
\end{align}
where $n>1$ is minimal. Note that the minimality of $n$ implies that the $c_i$'s and $p_i$'s are linearly independent over $K$.
\end{setup}


\subsection{Variety of all $K$-algebras}

The homogeneous variety $\mathcal{A}_K$ of all $K$-algebras  has the empty set of defining relations. Consider $K\{Y\}$, the free $K$-algebra on a set $Y$.
It is well known that every nonassociative word $w$ of degree at least two has unique representation in the form of a product of two nonassociative words. Hence one can introduce
a total ordering in the set $V(Y)$ of nonassociative words as follows. Order the words of length one (variables)
arbitrarily. Assuming that the words of length $n$, $n\geq 1$, have been already ordered in such a
way that words of smaller length  precede words of greater length, then given two words $w_1,w_2\in V(Y)$ of length
$n+1$ and represented as a product of two nonassociative words of lesser length $w_1=u_1v_1$, $w_2=u_2v_2$ we
define $w_1<w_2$ if, and only if, $u_1<u_2$ or $u_1=u_2$ and $v_1<v_2$.
Observe that this total ordering
of the non-associative words satisfies that if $w_1<w_2$ then $ww_1<ww_2$ for all $w\in V(Y)$.

\begin{theorem}
The variety $\mathcal{A}_K$ is MLM.
\end{theorem}

\begin{proof}
Consider $\mathfrak{A}=\mathcal{A}_K$ and $\mathfrak{A}_F=\mathcal{A}_F$ in Setup \ref{setup}.

Suppose $p_n$ is the $p_i$ of greatest degree $\geq 1$ and with greatest word in the support among the $p_i$'s.

If $\mu\left(\sum\limits_{i=1}^{n}{c_i\otimes p_i}\right)=\sum\limits_{i=1}^n c_ip_i =0,$
then,
\begin{eqnarray*}
0 & = &  \left(\sum_{i=1}^nc_ip_i\right) p_n-p_n\left(\sum_{i=1}^n c_ip_i \right) \\
& = &
\sum_{i=1}^n c_i\left(p_i p_n - p_n p_i\right) = \sum_{i=1}^{n-1} c_i\left(p_i p_n - p_n p_i\right)
\end{eqnarray*}

Hence,
\begin{equation*}
    \mu\left(\sum_{i=1}^{n-1} c_i\otimes \Big(p_i p_n - p_n p_i\Big)\right)=0.
\end{equation*}
By the minimality of $n$, we have
$0=p_i p_n-p_n p_i \textrm{ for }i=1,\dotsc,n.$
This equality implies that $p_i$ and $p_n$ have the same maximal word in its support. Thus, there exists $\lambda_i\in K$
such that $p_i-\lambda_i p_n$ has a lesser maximal word in its support. Since
\begin{equation*}
    (p_i-\lambda_i p_n)p_n -p_n(p_i-\lambda_i p_n)=0 \textrm{ for }i=1,\dotsc,n,
\end{equation*}
we obtain that $p_i=\lambda_ip_n$ for each $i=1,\dotsc,n.$ This contradicts the fact that the
$p_i$'s are $K$-linear independent.
\end{proof}


\subsection{Variety of commutative $K$-algebras}

The homogeneous variety $\mathcal{C}_K$ of commutative $K$-algebras has the 
defining relation $x_1x_2-x_2x_1$. The free commutative $K$-algebra on a set $Y$ will be denoted by $\mathcal{C}_K\{Y\}$.
We will need a result and a definition  from \cite{Shirshov}. 

Let $Y$ be a nonempty set. Consider the words on $Y$. Words of length one will be called
\emph{regular} and ordered arbitrarily. Assuming that regular words of length less than
$n$, $n>1$, have been already defined and ordered in such a way that words of smaller
length precede words of greater length, a word $w$ of length $n$ will be called \emph{regular} if
\begin{enumerate}
	\item $w=uv$ where $u$ and $v$ are regular words;
	\item $u\geq v$.
\end{enumerate}
We order the regular words of length $n$ defined in this way,  declaring that
$w_1=u_1v_1<w_2=u_2v_2$ if either $u_1<u_2$ or $u_1=u_2$ and $v_1<v_2$. Then
we declare the regular words of length $n$ to be greater
than regular words of smaller length.
By \cite[Theorem~1]{Shirshov}, the collection of all regular words form a basis
of  $\mathcal{C}_K\{Y\}$.

\begin{theorem}
The variety $\mathcal{C}_K$ is MLM.
\end{theorem}

\begin{proof}
Consider $\mathfrak{A}=\mathcal{C}_K$ and $\mathfrak{A}_F=\mathcal{C}_F$ in Setup \ref{setup}.

Suppose $p_n$ is the $p_i$ of greatest degree $\geq 1$ and with greatest word in the support among the $p_i$'s.

If $\mu\left(\sum\limits_{i=1}^{n}{c_i\otimes p_i}\right)=\sum\limits_{i=1}^n c_ip_i =0,$
then, for all regular words $w$,
\begin{eqnarray*}
0 & = &  \left(\left(\sum_{i=1}^nc_ip_i\right)w\right) p_n-(p_nw)\left(\sum_{i=1}^n c_ip_i \right) \\
& = &
\sum_{i=1}^n c_i\Big((p_iw) p_n - (p_nw) p_i\Big) =  \sum_{i=1}^{n-1} c_i\Big((p_iw) p_n - (p_nw). p_i\Big)
\end{eqnarray*}

Hence,
\begin{equation*}
    \mu\left(\sum_{i=1}^{n-1} c_i\otimes \Big((p_iw) p_n - (p_nw) p_i\Big)\right)=0,
\end{equation*}
for all regular words $w$.
By the minimality of $n$,
\begin{equation*}
    0=(p_iw) p_n-(p_nw) p_i \textrm{ for all regular words $w$ and for }i=1,\dotsc,n.
\end{equation*}
From the commutativity of the product and the definition of regular words,
this equality implies that $p_i$ and $p_n$ have the same maximal regular word in its support.
Thus, there exists $\lambda_i\in K$
such that $p_i-\lambda_i p_n$ has a lesser maximal regular word in its support. Since
$((p_i-\lambda_i p_n)w)p_n -(p_nw)(p_i-\lambda_i p_n)=0$ for all regular words $w$ and $i=1,\dotsc,n$,
we obtain that $p_i=\lambda_ip_n$ for each $i=1,\dotsc,n,$ a contradiction.
\end{proof}


\subsection{Schreier varieties satisfying $x^2=0$}

Let $\mathfrak{M}$ be a variety of $K$-algebras. We say that $\mathfrak{M}$
is a \emph{Schreier variety} if any subalgebra of a free algebra in $\mathfrak{M}$ is a free algebra in $\mathfrak{M}$.
The main
examples of homogeneous Schreier varieties are: the varieties of all algebras \cite{Kurosh47}, commutative and anticommutative algebras \cite{Shirshov}, Lie algebras \cite{Shirshov53} and \cite{Witt},
and algebras with zero multiplication.

Let $A$ be a free $K$-algebra in $\mathfrak{M}$. Let $(x_1,\dotsc,x_n),\, (y_1,\dotsc,y_n)\in A^n$. Let
$V$ and $W$ be the $K$-subspaces of $A$ generated by $(x_1,\dotsc,x_n)$ and $(y_1,\dotsc,y_n)$, respectively.
We say that a transformation $\tau\colon(x_1,\dotsc,x_n)\mapsto (\tau(x_1)=y_1,\dotsc,\tau(x_n)=y_n)$ is an
\emph{elementary transformation} if either:
\begin{enumerate}
	\item $\tau$ induces a nonsingular $K$-linear transformation between $V$ and $W$, or
	\item $y_1=x_1,y_2=x_2,\dotsc,y_{n-1}=x_{n-1}$ and $y_n=x_n+u_n$ where $u_n$ belongs to the $K$-subalgebra
	of $A$ generated by $x_1,\dotsc,x_{n-1}$.
\end{enumerate}
Observe that inverses of elementary transformations are also elementary transformations.

It is known that homogeneous Schreier varieties $\mathfrak{M}$ are \emph{Nielsen},
see for example \cite{Lewin} or  \cite[Chapter~11]{MiShYu}. In other words,
if $A$ is a free algebra of $\mathfrak{M}$,
one can transform any finite set of elements $a_1,\dotsc,a_n\in A$
to a free set of generators of the free subalgebra generated by $a_1,\dotsc,a_n$
by using a finite number of elementary transformations and cancelling possible zero elements.

\begin{lemma}\label{L1}
Suppose that $\mathfrak{M}$ is a homogeneous Schreier variety of $K$-algebras
that satisfies the identity $x^2=0$ and does not satisfies
the identity $x_1x_2=0$. Let $Y$ be a set of at least two elements and
denote by $K_\mathfrak{M}\{Y\}$ the free algebra in the variety $\mathfrak{M}$ with set
of free generators $Y$. Suppose that $p,q\in K_\mathfrak{M}\{Y\}$, $q\neq 0$, such that $pq=0.$ Then
$p=\lambda q$ for some $\lambda\in K$.
\end{lemma}

\begin{proof}
Since $\mathfrak{M}$ is a Schreier variety, the subalgebra $B$ of $K_\mathfrak{M}\{Y\}$
generated by $p,q$ is a free subalgebra of $K_\mathfrak{M}\{Y\}$. There exist
elementary transformations
\begin{equation*}
    (p,q)\rightarrow \dotsb \rightarrow (x,y),
\end{equation*}
such that $x,y$ (and $y$ may be zero) are free generators of $B$.
If the rank of $B$ is two, then there exist elementary transformations (the inverses of the previous ones)
\begin{equation*}
    (x,y)\rightarrow \dotsb \rightarrow (p,q).
\end{equation*}

Since $B$ is free on $x,y$, the elementary transformations induce automorphisms of $K$-algebras.
Hence there exists an isomorphism of $B$ that sends $x\mapsto p$, $y\mapsto q$. Hence
$p,q$ are free generators of $B$. Since $pq=0$, and $p^2=q^2=0$, it implies that
the free algebra in the variety $\mathfrak{M}$ on two free generators has zero product.
Therefore, all the algebras in $\mathfrak{M}$ satisfy the identity $x_1x_2=0$, a contradiction.
It implies that $B$ is a free algebra of rank one. Since $x^2=0$, the free subalgebra generated
by one nonzero element is of dimension $1$. Now $B$ is generated by $q$, as desired.
\end{proof}

  Between the homogeneous Schreier varieties that satisfy the conditions of Lemma \ref{L1}, we highlight the following:
\begin{itemize}
    \item The variety of Lie $K$-algebras  $\mathcal{L}_K$. Its defining relations are: $x^2_1$ and $(x_1x_2)x_3+(x_2x_3)x_1+(x_3x_1)x_2$. 
    \item The variety of anticommutative $K$-algebras $\mathcal{AC}_K$. It has the defining relation: $x_1x_2+x_2x_1$. 
\end{itemize}

\begin{theorem}\label{theo:Schreier}
Let $\mathfrak{M}$ be a Schreier variety of $K$-algebras
  that satisfies the identity $x^2=0$ and does not satisfy
	the identity $x_1x_2=0$. Then $\mathfrak{M}$ is MLM. In particular,
the varieties $\mathcal{L}_K$ and $\mathcal{AC}_K$ are MLM.
\end{theorem}

\begin{proof}
Consider $\mathfrak{A}=\mathfrak{M}$ and $\mathfrak{A}_F=\mathfrak{M}_F$ in Setup \ref{setup}. Then
\begin{align*}
\mu\left(\sum_{i=1}^{n-1} c_i\otimes p_ip_n\right)=\mu\left(\sum_{i=1}^{n} c_i\otimes p_ip_n\right)=
\sum_{i=1}^n c_i p_ip_n=\left(\sum_{i=1}^{n}c_ip_i \right)p_n=0\cdot p_n=0.
\end{align*}
By the minimality of $n$, $p_ip_n=0$ for all $i=1,\dotsc,n$. By Lemma~\ref{L1}, and the fact that $p_n\neq 0$,
we get that $p_i=\lambda_i p_n$, where $\lambda_i\in K$ for $i=1,\dotsc,n$, a contradiction.
\end{proof}


\subsection{Variety generated by special Jordan algebras}

Let $K$ be a field.
If $A$ is an associative $K$-algebra,  we define on the $K$-vector space $A$ a new multiplication $\circ$,
which is connected with the associative multiplication by the formula
$x\circ y=\frac{1}{2}(xy+yx)$.
In this way a new $K$-algebra is obtained and it is denoted by $A^{(+)}$. The $K$-algebra $A^{(+)}$
is a Jordan algebra (that is, satisfies $x_1x_2=x_2x_1$ and $(x_1^2x_2)x_1=x_1^2(x_2x_1)$).
If $J$ is a $K$-subspace of $A$ which is closed with respect to the operation
$\circ$, then $J$ together with $\circ$ is a subalgebra of $A^{(+)}$,
and consequently a Jordan algebra.  Such a Jordan algebra is called special Jordan algebra.
The variety generated by all special Jordan $K$-algebras will be denoted by $\mathcal{SJ}_K$ and it
consists of all  $K$-algebras which can be obtained as homomorphic
images of special Jordan $K$-algebras.

Let $Y$ be a set. Consider $K\langle Y\rangle$, the free associative $K$-algebra on the set $Y$. The
$K$-subalgebra  $\mathcal{SJ}_K(Y)$ of the algebra $K\langle Y\rangle^{(+)}$ generated by the set $Y$ is the free special Jordan $K$-algebra on $Y$ and it is the
free $K$-algebra on the set $Y$ in the variety $\mathcal{SJ}_K$.

\begin{theorem}
Let $K$ be a field of characteristic not two. Then  $\mathcal{SJ}_K$ is MLM.
\end{theorem}

\begin{proof}
Consider $\mathfrak{A}=\mathcal{SJ}_K$ and $\mathfrak{A}_F=\mathcal{SJ}_F$ in Setup \ref{setup}.

Suppose $p_n$ is the $p_i$ of greatest degree. If it is of degree zero, then all $p_i$ are
of degree zero and the result follows from the case $n=1$ because there exist $a_i\in K$ such that
\begin{align*}
\mu\left(\sum_{i=1}^{n}{c_i\otimes p_i}\right)=\mu\left(\sum_{i=1}^n c_ia_i\otimes 1\right)=0
\end{align*}

From now on, we suppose that $p_n$ is of positive degree.

\underline{Claim~1:} $p_i$ and $p_n$ commute as elements of $K\langle Y\rangle$.

If $\mu\left(\sum\limits_{i=1}^{n}{c_i\otimes p_i}\right)=\sum\limits_{i=1}^n c_ip_i =0,$
then, for all $z\in \mathcal{SJ}_K(Y)$
\begin{eqnarray*}
0 & = &  \left(\left(\sum_{i=1}^nc_ip_i\right)\circ z\right)\circ p_n-\left(\sum_{i=1}^n c_ip_i \right)\circ (z\circ p_n)  \\
& = &
\sum_{i=1}^n c_i\left((p_i\circ z)\circ p_n - (p_i\circ (z\circ p_n))\right) = \sum_{i=1}^{n-1} c_i\left((p_i\circ z)\circ p_n - (p_i\circ (z\circ p_n))\right).
\end{eqnarray*}

Hence,
\begin{equation*}
    \mu\left(\sum_{i=1}^{n-1} c_i\otimes \Big((p_i\circ z)\circ p_n - (p_i\circ (z\circ p_n))\Big)\right)=0.
\end{equation*}

By the minimality of $n$,
$0=(p_i\circ z)\circ p_n-p_i\circ (z\circ p_n)$ for $i=1,\dotsc,n$.
Using that $\mathcal{SJ}_K(Y)$ is a subalgebra of $K\langle Y\rangle ^{(+)}$, we get
\begin{eqnarray*}
0 & =& (p_i\circ z)\circ p_n-p_i\circ (z\circ p_n)\\ &=& \frac{1}{2}((zp_i+p_iz)\circ p_n-p_i\circ(zp_n+p_nz))\\
& =& \frac{1}{4}( zp_ip_n+p_izp_n+p_nzp_i+p_np_iz- (p_izp_n+p_ip_nz+zp_np_i+p_nzp_i)) \\
& = & \frac{1}{4}[z(p_ip_n-p_np_i)-(p_ip_n-p_np_i)z].
\end{eqnarray*}
This implies that $p_ip_n-p_np_i$ commutes with any variable $z\in Y\subseteq \mathcal{SJ}_K(Y)$.
Hence $p_ip_n-p_np_i$ is in the center of $K\langle Y\rangle $, that is, $p_ip_n-p_np_i\in K$. But
this is only possible if $p_ip_n-p_np_i=0$, because the term of zero degree in $p_ip_n-p_np_i$ is zero.
Therefore $p_i$ and $p_n$ commute in $K\langle Y\rangle$ and the claim is proved.

\underline{Claim~2:} There exists $\lambda_i \in K$ such that
  $p_i=\lambda_i p_n$ for $i=1,\dotsc,n$. A contradiction with the fact that the $p_i$'s are $K$-linear independent.

By \cite[Corollary~6.7.7]{CohnFreeIdeal} and Claim~1, there exists a polynomial $u\in K\langle Y\rangle$
of degree at least one such that $p_i\in k[u]$, $i=1,\dotsc,n$.

From \eqref{eq:muzero}, $\sum\limits_{i=1}^n c_ip_i=0$. Thus, for all $z\in\mathcal{SJ}_K(Y)$,
\begin{eqnarray*}
0& = &\left(\sum_{i=1}^n c_ip_i\right)\circ (z\circ p_n^2)-\left(\left(\sum_{i=1}^n c_ip_i\right)\circ p_n \right)\circ(z\circ p_n)\\
 & = & \sum_{i=1}^n c_i\Big( (p_i\circ (z\circ p_n^2))-(p_i\circ p_n)\circ (z\circ p_n)  \Big) \\
& = & \sum_{i=1}^{n-1} c_i \Big( (p_i\circ (z\circ p_n^2))-(p_i\circ p_n)\circ (z\circ p_n)  \Big),
\end{eqnarray*}
where we have used that $(z\circ p_n^2)\circ p_n=p_n^2\circ (z\circ p_n)$  in the last equality.

Hence
\begin{align*}
\mu\left( \sum_{i=1}^{n-1} c_i\otimes \left[(p_i\circ (z\circ p_n^2))-(p_i\circ p_n)\circ (z\circ p_n)\right]\right)=0.
\end{align*}
By the minimality of $n$, and using Step~1 and $\mathcal{SJ}_K(Y)\subseteq K\langle Y\rangle ^{(+)}$,
\begin{eqnarray*}
0 & = & p_i\circ (z\circ p_n^2) - (p_i\circ p_n) \circ (z\circ p_n) \\
& = & \frac{1}{2}p_i\circ (zp_n^2+p_n^2z)-\frac{1}{4}(p_ip_n+p_np_i)\circ (zp_n+p_nz) \\
& = & \frac{1}{4}(p_izp_n^2+p_ip_n^2z+zp_n^2p_i+p_n^2zp_i)\\
&  & -\frac{1}{4}(p_ip_nzp_n+p_ip_n^2z+zp_np_ip_n+p_nzp_ip_n) \\
& = & \frac{1}{4}(p_izp_n^2+p_n^2zp_i-p_ip_nzp_n-p_nzp_ip_n)\\
& = & \frac{1}{4}(p_i(zp_n-p_nz)p_n-p_n(p_nz-zp_n)p_i).
\end{eqnarray*}
Hence, $p_i(zp_n-p_nz)p_n=p_n(zp_n-p_nz)p_i$ for all  $z\in \mathcal{SJ}_K(Y)\subseteq K\langle Y\rangle ^{(+)}$ and $i=1,\dotsc,n$. 

Let now $x\in \mathcal{SJ}_K(Y)$ such that $x$ does not commute with $u$. By \cite[Corollary~6.7.4]{CohnFreeIdeal}, $u$ and $x$ form a free set over $K$.
Then, rewriting the last equality, we obtain
    $p_i(u)(xp_n(u)-p_n(u)x)p_n(u)=p_n(u)(xp_n(u)-p_n(u)x)p_i(u)$.

Suppose that the degree on $u$ of $p_i(u)$ is smaller than the degree of $p_n(u)$. Considering the lexicographic order
in $K\langle x,u\rangle$ with $x<u$, we obtain that the greatest monomial on the right hand side of the equation
is obtained from $-p_n(u)^2xp_i(u)$, which cannot be obtained from any other monomial on the left hand side of the equality. A contradiction.

Hence $p_i$ and $p_n$ are polynomials on $u$ of the same degree. Thus there exists $\lambda\in K$ such that
$p_i(u)-\lambda p_n(u)$ has degree smaller than the degree of $p_n$ and
\begin{equation*}
    (p_i(u)-\lambda p_n(u))(xp_n(u)-p_n(u)x)p_n(u)=p_n(u)(xp_n(u)-p_n(u)x)(p_i(u)-\lambda p_n(u)).
\end{equation*}
By the foregoing, the only possibility is that $p_i=\lambda_i p_n$ for some $\lambda_i\in K$.
\end{proof}

\subsection{Variety of non-commutative Poisson $K$-algebras}
In this subsection we deviate from our notation and consider a variety of algebras endowed with more than one product. More precisely,  the variety $\mathcal{NCP}_K$ of all non-commutative Poisson $K$-algebras.

A vector space $P$ over a field $K$ endowed with two bilinear operations $x\cdot y$ (a multiplication) and $\{x, y\}$ (a Poisson
bracket) is called a non-commutative Poisson algebra if $P$ is an associative algebra under $x\cdot y$, $P$ is a Lie algebra under $\{x, y\}$, and $P$ satisfies the  Leibniz identity:
$\{x \cdot y, z\} = \{x, z\} \cdot y + x \cdot \{y, z\}$ for all $x,y,z\in P$.

Let $Y$ be a set. 
The free non-commutative Poisson $K$-algebra $\mathcal{P}_K(Y)$ is constructed as follows. Let $\mathcal{L}_K(Y)$ be the free Lie $K$-algebra on $Y$ and suppose that
$X=\{x_1,x_2,\dotsc\}$ is a $K$-basis of $\mathcal{L}_K(Y)$. Then $\mathcal{P}_K(Y)$ is the
free associative algebra on the set of free generators $X$.  Using the Leibniz identity
one can uniquely extend the Lie bracket $\{x, y\}$ of $\mathcal{L}_K(Y)$ to a Poisson bracket $\{x, y\}$ on
$\mathcal{P}_K(Y)$, and $\mathcal{P}_K(Y)$ becomes a Poisson algebra.

\begin{corollary}
The variety $\mathcal{NCP}_K$ is MLM.
\end{corollary}
\begin{proof}
Let $K\subseteq F$ be a field extension and let $A\in\mathcal{NCP}_F$. Suppose that $A$ contains a free non-commutative Poisson $K$-algebra on a free set of generators $Y\subseteq A$ of at least two elements. The free Lie $K$-algebra (with respect to $\{x,y\}$) generated by $Y$ is the free Lie $K$-algebra $\mathcal{L}_K(Y)$. By Theorem~\ref{theo:Schreier}, the Lie $F$-subalgebra of $A$ generated by $Y$ is the free Lie $F$-algebra $\mathcal{L}_F(Y)$ with  set of free generators $Y$. Note that we can pick the same basis $\mathcal{B}$ for $\mathcal{L}_K(Y)$ and
$\mathcal{L}_F(Y)$.
Now, by Lemma~\ref{lem:MLM}, the associative $F$-subalgebra generated by $\mathcal{B}$ is the free associative $F$-algebra on $\mathcal{B}$, as desired.
\end{proof}

We would like to finish noting that one can mimic the proof of Theorem~\ref{theo:MLMReichstein} to show that the variety 
$\mathcal{NCP}_K$ is Reichstein when $K$ is an uncountable field.


\end{document}